\documentclass[a4paper, 12pt, leqno]{amsart}
\usepackage{amsmath}
\usepackage{amscd}
\usepackage{amssymb}
\usepackage{amsthm}
\newtheorem{lemma}{{\sc Lemma}}[section]

\newtheorem{proposition}[lemma]{{\sc Proposition}}
\newtheorem{theorem}[lemma]{{\sc Theorem}}
\newtheorem{remark}[lemma]{{\sc Remark}}
\newtheorem{conjecture}[lemma]{{\sc Conjecture}}

\numberwithin{equation}{section}

\def\Gb{{\mathfrak{b}}}
\def\Gg{{\mathfrak{g}}}

\def\Gn{{\mathfrak{n}}}
\def\Gt{{\mathfrak{t}}}
\def\Gz{{\mathfrak{z}}}
\def\Gsl{{\mathfrak{sl}}}
\def\Ggl{{\mathfrak{gl}}}
\def\BC{{\mathbb{C}}}
\def\BF{{\mathbb{F}}}

\def\BLL{{\mathbb{L}}}
\def\BQ{{\mathbb{Q}}}
\def\BZ{{\mathbb{Z}}}

\def\CB{{\mathcal B}}
\def\CC{{\mathcal C}}
\def\CO{{\mathcal O}}
\def\CI{{\mathcal I}}
\def\CL{{\mathcal L}}
\def\CN{{\mathcal N}}
\def\CP{{\mathcal P}}
\def\CQ{{\mathcal Q}}
\def\CR{{\mathcal R}}

\def\Ad{\mathop{\rm Ad}\nolimits}

\def\Aut{{\rm Aut}}
\def\Coh{\mathop{\rm Coh}\nolimits}
\def\Coker{\mathop{\rm Coker}\nolimits}

\def\gr{\mathop{\rm{gr}}\nolimits}

\def\poi{{\mathop{\rm pt}\nolimits}}
\def\Supp{{\rm{Supp}}}
\def\Tor{{\rm{Tor}}}
\def\ll{\underset L\leq}
\def\rl{\underset R\leq}
\def\lr{\rl}
\def\lrl{\underset {LR}\leq}
\def\llr{\lrl}
\def\el{\underset L\sim}
\def\er{\underset R\sim}
\def\elr{\underset {LR}\sim}
\def\ga{\gamma}

\title[Kazhdan-Lusztig Basis  and A Geometric Filtration
]{Kazhdan-Lusztig Basis and A Geometric Filtration of  an Affine
Hecke Algebra}
\date{\today}
\author{Toshiyuki TANISAKI$^{*}$ and Nanhua XI$^{\dagger}$}
\address{$^{*}$
Department of Mathematics\\
Osaka City University \\
3-3-138, Sugimoto, Sumiyoshi-ku \\
Osaka, 558-8585\\
Japan} \email{tanisaki@sci.osaka-cu.ac.jp}
\address{$^{\dagger}$
Institute of Mathematics\\
Chinese Academy of Sciences\\
Beijing, 100080\\
China } \email{nanhua@math.ac.cn}
\thanks{N. Xi was partially supported by a fund of 973
program.}

\begin{document}
\begin{abstract}
According to Kazhdan-Lusztig and Ginzburg, the Hecke algebra of an
affine Weyl group is identified with the equivariant $K$-group of
Steinberg's triple variety. The $K$-group is equipped with a
filtration indexed by closed $G$-stable subvarieties of the
nilpotent variety, where $G$ is the corresponding reductive
algebraic group over $\BC$. In this paper we will show in the case
of type $A$ that the filtration is compatible with the
Kazhdan-Lusztig basis of the Hecke algebra.
\end{abstract}
\maketitle \setcounter{section}{-1}
\section{Introduction}
\label{sec:Intro} Let $G$ be a connected reductive algebraic group
over the complex number field $\BC$ with simply-connected derived
group. Let $W$ and $P$ be its Weyl group and weight lattice
respectively. The semidirect product $\tilde{W}_a=WP$ with respect
to the action of $W$ on $P$ is called an (extended) affine Weyl
group. Let $H(\tilde{W}_a)$ be the associated Hecke algebra.
According to Kazhdan-Lusztig and Ginzburg (\cite{KL2, G}) we have
a geometric realization of $H(\tilde{W}_a)$ in terms of
equivariant $K$-theory. Namely, we have an isomorphism
\[
\Phi:H(\tilde{W}_a)\to K^{G\times\BC^*}(Z)
\]
of $\BZ[q^{1/2},q^{-1/2}]$-algebras, where $K^{G\times\BC^*}(Z)$
denotes the equivariant $K$-group of Steinberg's triple variety
$Z$ with respect to the natural action of $G\times\BC^*$. Let
$\CN$ be the nilpotent variety of the Lie algebra $\Gg$ of $G$.
For each $G$-stable closed subset $V$ of $\CN$ there corresponds a
$G\times\BC^*$-stable closed subvariety $Z_V$ of $Z$, and the
associated equivariant $K$-group $K^{G\times\BC^*}(Z_V)$ is
identified with a two-sided ideal of $K^{G\times\BC^*}(Z)$.
Moreover, we have $K^{G\times\BC^*}(Z_{V_1})\subset
K^{G\times\BC^*}(Z_{V_2})$ if $V_1\subset V_2$.

Recall  that $H(\tilde{W}_a)$ is equipped with the Kazhdan-Lusztig
basis  $\{C_w\mid w\in \tilde W_a\}$ (\cite{KL1}). It plays very
important roles in various aspects of the representation theory of
reductive algebraic groups. It should be an interesting problem to
give a geometric description of $\Phi(C_w)$ for $w\in\tilde{W}_a$.
An answer in the case $w\in W$ is given in \cite{T}. Moreover, the
answer for certain elements corresponding to  dominant elements in
$P$ is given in \cite{L1}. Related to this problem, it is
conjectured that $K^{G\times\BC^*}(Z_V)$ is spanned by a subset of
$\{\Phi(C_w)\mid w\in \tilde W_a\}$ for any $G$-stable closed
subset $V$ of $\CN$. In particular, any $H(\tilde{W}_a)$-bimodule
associated to a two-sided cell of $\tilde{W}_a$ should be
identified with $ K^{G\times\BC^*}(Z_{\overline{O}})/
K^{G\times\BC^*}(Z_{\overline{O}\setminus O}) $ for a nilpotent
orbit $O$.

The aim of this paper is to prove this conjecture in the case $G$
is of type $A$. A key to this result is the fact that the
$H(\tilde{W}_a)$-bimodule corresponding to a two-sided cell of
$\tilde{W}_a$ is generated by a single element (see Theorem
\ref{thm:key1} below).

The contents of this paper are as follows. In Section
\ref{sec:Hecke} and Section \ref{sec:affine Hecke} we will recall
some fundamental facts on (affine) Hecke algebras. A precise
formulation of the above stated conjecture in view of the
bijection between the set of nilpotent orbits and that of
two-sided cells will be given in Section \ref{sec:HK}. In Section
\ref{sec:GL} we will give a proof of the conjecture in the case
$G=GL_n(\BC)$. The arguments works for $SL_n(\BC)$ as well. In
Appendix A we will collect well-known facts on equivariant
$K$-theory, and in Appendix B we will give a description of the
product on the quotient $ K^{G\times\BC^*}(Z_{\overline{O}})/
K^{G\times\BC^*}(Z_{\overline{O}\setminus O})$ for any $G$ in
terms of the Springer fiber and Slodowy's variety, where $O$ is a
nilpotent orbit.

\section{Hecke algebras}
\label{sec:Hecke} Let $(W,S)$ be a Coxeter system with the length
function $\ell:W\to\BZ_{\geq0}$ and the standard partial order
$\geq$. Assume that we are given a group $\Omega$ and a group
homomorphism $\Omega\to\Aut(W,S)$, where $\Aut(W,S)$ denotes the
automorphism group of $(W,S)$. We denote by $\tilde{W}$ the
semidirect product $\Omega W$ with respect to the action of
$\Omega$ on $W$. The length function $\ell$ and the standard
partial order $\geq$ for $W$ are naturally extended to $\tilde{W}$
by
\begin{align*}
&\ell(\omega w)=\ell(w) \qquad
(\omega\in\Omega, w\in W),\\
&\omega_1w_1\geq\omega_2w_2 \Leftrightarrow
\omega_1=\omega_2,\,\,w_1\geq w_2
\qquad(\omega_1,\omega_2\in\Omega, w_1, w_2\in W).
\end{align*}
For $w$ in $\tilde{W}$ we set
\[
L(w)=\{s\in S\ |\ sw\leq w\},\qquad R(w)=\{s\in S\ |\ ws\leq w\}.
\]

We denote by $H(\tilde{W})$ the Hecke algebra associated to
$\tilde{W}$ . It is an associative algebra over the Laurent
polynomial ring $\BZ[q^{1/2},q^{-1/2}]$. As a
$\BZ[q^{1/2},q^{-1/2}]$-module it has a free basis $\{T_w\mid w\in
\tilde{W}\}$, and the multiplication is determined by
\begin{align*}
&T_yT_w=T_{yw}\qquad
(y, w\in \tilde{W}, \,\,\ell(y)+\ell(w)=\ell(yw)),\\
&(T_s+1)(T_s-q)=0\qquad (s\in S).
\end{align*}
There is a unique ring automorphism $h\mapsto\overline{h}$ of
$H(\tilde{W})$  determined by
\[
\overline{q^{1/2}}=q^{-1/2},\qquad
\overline{T}_w=T_{w^{-1}}^{-1}\quad(w\in \tilde{W})
\]

According to Kazhdan-Lusztig \cite{KL1}, for each $w\in \tilde{W}$
there exists uniquely an element
\[
C_w=\sum_{y\leq w}P_{y,w}(q)T_y
\]
of $H(\tilde{W})$ satisfying
\begin{itemize}
\item[(a)] $P_{w,w}(q)=1$, \item[(b)] for $y<w$ we have
$P_{y,w}(q)\in\BZ[q]$, and
$\deg(P_{y,w}(q))\leq(\ell(w)-\ell(y)-1)/2$, \item[(c)]
$\overline{C}_w=q^{-\ell(w)}C_w$.
\end{itemize}
The basis $\{C_w\mid w\in \tilde{W}\}$ of $H(\tilde{W})$ is called
the Kazhdan-Lusztig basis. We will also use
\[
C'_w=q^{-{\ell(w)}/2}C_w\qquad(w\in\tilde{W}).
\]

For $w\in\tilde{W}$ let $\CI_w$ (resp. $\CI_w^L$, $\CI_w^R$)
denote the set of two-sided (resp. left, right) ideals $I$ of
$H(\tilde{W})$ subject to the conditions
\begin{itemize}
\item[(a)] $C_w\in I$, \item[(b)] $I$ is spanned over
$\BZ[q^{1/2},q^{-1/2}]$ by a subset of $\{C_y\mid
y\in\tilde{W}\}$.
\end{itemize}
It contains the unique minimal element $I_w=\bigcap_{I\in \CI_w}I$
(resp. $I_w^L=\bigcap_{I\in \CI_w^L}I$, $I_w^R=\bigcap_{I\in
\CI_w^R}I$). We define a preorder $\lrl$ (resp. $\ll, \ \rl$) and
an equivalence relation $\elr$ (resp. $\el,\ \er$) on $\tilde{W}$
by
\begin{align*}
&y\lrl w\Leftrightarrow I_y\subset I_w,
\\
&{\text{(resp. }} y\ll w \Leftrightarrow I^L_y\subset I^L_w, \quad
y\rl w \Leftrightarrow I^R_y\subset I^R_w)
\\ &y\elr
w\Leftrightarrow I_y= I_w,\\
  &{\text{(resp. }} y\el w \Leftrightarrow I^L_y= I^L_w, \quad
y\er w \Leftrightarrow I^R_y= I^R_w).
\end{align*}

Equivalence classes with respect to $\elr $ (resp. $\el,\ \er$)
are called two-sided (resp. left, right) cells of $\tilde{W}$. The
preorder $\lrl $ on $\tilde{W}$ induces a partial order on the set
of two-sided cells which is also denoted by $\lrl $. For a
two-sided cell $\CC$ of $\tilde{W}$ with $w\in\CC$ we define 
two-sided ideals $H(\tilde{W})_{\lrl \CC}$ and
$H(\tilde{W})_{\underset{LR}<\CC}$ of $H(\tilde W)$ by
\[
H(\tilde{W})_{\lrl  \CC}=I_w,\qquad
H(\tilde{W})_{\underset{LR}<\CC}=\sum_{y\lrl w,\ y\notin\CC}I_y.
\]
The $H(\tilde{W})$-bimodule
\[
H(\tilde{W})_\CC=H(\tilde{W})_{\lrl  \CC}/
H(\tilde{W})_{\underset{LR}<\CC}\] has a canonical
$\BZ[q^{1/2},q^{-1/2}]$-basis parametrized by $\CC$. The
multiplication of $H(\tilde{W})$ induces a multiplication of
$H(\tilde{W})_\CC$ which is associative; however,
$H(\tilde{W})_\CC$ does not contain the identity element in
general.

\begin{lemma}
[Kazhdan-Lusztig \cite{KL1}] \label{lem:LR} If $y\ll w$ $($resp.
$y\lr w$$)$, then $R(w)\subset R(y)$ $($resp. $L(w)\subset
L(y)$$)$. In particular, if  $y\el w$ $($resp. $y\er w$$)$, then
$R(w)=R(y)$ $($resp. $L(w)=L(y)$$)$.
\end{lemma}

For a subset $T$ of $S$ such that $\langle T\rangle$ is a finite
subgroup of $W$ we denote the longest element of $\langle
T\rangle$ by $w_T$. We call $w\in \tilde W$ a parabolic element if
there exists some $T\subset S$ such that $|\langle
T\rangle|<\infty$ and $w=w_T$.

We will need the following simple assertion later.
\begin{lemma}
\label{lem:pa} Let $x,y\in \tilde W$ and let $w$ be a parabolic
element of $\tilde W$. Assume that  $x\ll w$ and $y\rl w$. Then
$C'_x=hC'_w$ and $C'_y=C'_wh'$ for some $h,h'\in H(\tilde W)$.
\end{lemma}
\begin{proof}
By Lemma \ref{lem:LR} we have $R(w)\subset R(x)$ and $L(w)\subset
L(y)$. Since $w$ is a parabolic element, there are $x_1$ and $y_1$
in $\tilde W$ such that $x=x_1w$, $y=wy_1$ and
$\ell(x)=\ell(x_1)+\ell(w)$, $\ell(y)=\ell(w)+\ell(y_1)$.
  Now
using induction on the length of $x_1$ and of $y_1$ we see the
assertion is true (see \cite[(2.3.a), (2.3.b)]{KL1}).
\end{proof}
In the analysis of two-sided cells the star operations defined in
Kazhdan-Lusztig \cite{KL1} and the $a$-function defined in Lusztig
\cite{L:cell1} play important roles.

Let $s$ and $t$ be in $S$ such that $st$
  has order 3, i.e. $sts=tst$. Define
$$D_L(s,t)=\{w\in \tilde W\  |\   L(w)\cap\{s,t\}\text{ has exactly one
element}\},$$
$$D_R(s,t)=\{w\in \tilde W\ |\ R(w)\cap\{s,t\}\text{ has
exactly one element}\}.$$ If $w$ is in $D_L(s,t)$, then
$\{sw,tw\}$ contains exactly one element in $D_L(s,t)$, denoted by
${}^*w$, here $*=\{s,t\}$. The map
\[
D_L(s,t)\ni w\mapsto{}^*w\in D_L(s,t)
\]
is an involution and is called a left star operation. Similarly we
can define the right star operation $D_R(s,t)\ni w\mapsto w^*\in
D_R(s,t)$ by $\{w^*\}=\{ws,wt\}\cap D_R(s,t)$.
\begin{proposition}
[Kazhdan-Lusztig \cite{KL1}] \label{prop:LRstar} Let $s$ and $t$
be in $S$ such that $st$
  has order 3, and set $*=\{s, t\}$.
\begin{itemize}
\item[(i)] For $w\in D_L(s,t)$ $($resp. $D_R(s,t)$$)$ we have
${}^*w\el w$ $($resp. $w^*\er w$$)$. \item[(ii)] For $y,w\in
D_L(s,t)$ $($resp. $D_R(s,t)$$)$ with $y\er w$ $($resp. $y \el
w$$)$ we have ${}^*y\er{}^*w$ $($resp. $y^*\el w^*$$)$.
\end{itemize}
\end{proposition}

Given $w,u$ in $\tilde W$, we write
$$C'_wC'_u=\sum_{v\in \tilde W_a}h_{w,u,v}C'_v\qquad(h_{w,u,v}\in
\BZ[q^{1/2},q^{-1/2}]).$$ The $a$-function
\[
a:\ \tilde W\to\BZ_{\geq0}\sqcup\{\infty\}
\]
is defined as follows. Let $v\in\tilde W$. If for any
$i\in\BZ_{\geq0}$ there exist some $w,u\in \tilde W$ such that
$q^{-i/2}h_{w,u,v}\not\in\BZ[q^{-1/2}]$, then we set
$a(v)=\infty$. Otherwise we set
  \[
  a(v)={\rm min}\{i\in\BZ_{\geq0}\ |\ q^{-
i/2}h_{w,u,v}\in\BZ[q^{-
  1/2}]{\rm\ for \ all\ }w,u\in \tilde W\}.
\]
\begin{proposition}
[Lusztig \cite{L:cell1}] \label{prop:a-crys} Assume that $(W,S)$
is crystallographic. Then for any $w,y\in \tilde W$ with $y\llr w$
we have $a(y)\geq a(w)$. In particular, the function $a$ is
constant on each two-sided cell of $\tilde W$

\end{proposition}

\section{Affine Hecke algebras}
\label{sec:affine Hecke} Let $G$ be a connected reductive
algebraic group over $\BC$ with simply-connected derived group.
Let $B$ and $T$ be a Borel subgroup and a maximal torus of $G$
respectively such that $B\supset T$. We denote the Lie algebras of
$G, B, T$ by $\Gg, \Gb, \Gt$ respectively. Let
$\Delta\subset\Gt^*$ be the root system. For $\alpha\in\Delta$ we
denote the corresponding root subspace by $\Gg_\alpha$. We choose
a system $\Delta^+$ of positive roots as the weights of $\Gg/\Gb$,
and denote the corresponding set of simple roots by $\Pi$. Let
$P\subset\Gt^*$ denote the weight lattice and let $Q$ be its
sublattice spanned by $\Delta$. We denote the subset of $P$
consisting of dominant weights by $P^+$.

In the rest of this paper we denote by $W$ the Weyl group of $G$.
It is a Coxeter group with canonical generator system
$S=\{s_{\alpha}\mid \alpha\in \Pi\}$. Here, the reflection with
respect to $\alpha\in\Delta$ is denoted by $s_\alpha$.

Let $W_a=WQ$ (resp. $\tilde{W}_a=WP$) denote the semidirect
product with respect to the action $W$ on $Q$ (resp. $P$). $W_a$
and $\tilde{W}_a$ are called the affine Weyl group and the
extended affine Weyl group respectively. The element of $W_a$
(resp. $\tilde{W}_a$ ) corresponding to $\lambda\in Q$ (resp.
$\lambda\in P$) is denoted by $t_\lambda$. Let $\Delta_c$ denote
the set of roots $\beta$ such that the corresponding coroots
$\beta^\vee$ are the highest coroots of irreducible components of
the coroot system, and set
\[
S_a=S\sqcup\{t_{\beta}s_{\beta} \mid\beta\in\Delta_c\}.
\]
Then $(W_a,S_a)$ is a Coxeter system. Set
\[
\Omega=\{\omega\in \tilde{W}_a\mid \omega S_a=S_a\omega\}.
\]
Then $\tilde{W}_a$ is canonically isomorphic to the semidirect
product $\Omega W_a$ with respect to the conjugation action of
$\Omega$ on $W_a$. Especially, we have the Hecke algebra
$H(\tilde{W}_a)$ of $\tilde{W}_a$. We identify the Hecke algebra
$H(W)$ of $W$ with a subalgebra of $H(\tilde W_a)$ by the
canonical embedding $T_w\mapsto T_w\,\,(w\in W)$.

We have the following properties on the $a$-function on $\tilde
W_a$.
\begin{proposition}
[Lusztig \cite{L:cell1}]
\begin{itemize}
\item[] \item[(i)] $a(w)=\ell(w)$ if $w$ is a parabolic element of
$\tilde W_a$. \item[(ii)] $a(w)\le a(w_S)(=\ell(w_S))$ for any
$w\in \tilde W_a$.
\end{itemize}
\end{proposition}
Note also that the function $a$ is constant on each two-sided cell
of $\tilde W_a$ by Proposition \ref{prop:a-crys}. For $w, u,
v\in\tilde{W}$ we define $\gamma_{w,u,v}\in\BZ$ by
$$h_{w,u,v}=\gamma_{w,u,v}q^{{a(v)}/2}+\text{lower degree
terms.}$$ Now we present a property of $\ga_{w,u,v}$ related to
  the star operations.

  By a similar argument to that for Theorem 1.4.5 in \cite{X:b}, we can see the following result.

\begin{proposition}
 \label{prop:star}
  Let $s,t$ be in $S_a$  such
that  $st$ has order 3. Set $*=\{s,t\}$. Assume $w, v\in
D_L(s,t)$. Then we have
\[
\ga_{w,u,v}=\ga_{{}^*w,u,{}^*v}.
\]
\end{proposition}

For $\lambda\in P$ we define $\theta_\lambda\in H(\tilde{W}_a)$ as
follows. Take $\lambda_1, \lambda_2\in P^+$ such that
$\lambda=\lambda_1-\lambda_2$ and set
\[
\theta_\lambda= q^{(-\ell(t_{\lambda_1})+\ell(t_{\lambda_2}))/2}
T_{t_{\lambda_1}} T_{t_{\lambda_2}}^{-1}.
\]
It does not depend on the choice of $\lambda_1, \lambda_2$.
Moreover, we have
\begin{align*}
&\theta_0=1,\\
&\theta_\lambda\theta_\mu=\theta_{\lambda+\mu}
\qquad(\lambda, \mu\in P),\\
&T_{s_\alpha}\theta_\lambda =\theta_{s\lambda}T_{s_\alpha}+(q-1)
\frac{\theta_\alpha(\theta_\lambda-\theta_{s\lambda})}
{\theta_\alpha-1} \qquad (\lambda\in P, \alpha\in\Pi).
\end{align*}
This presentation in terms of the $\BZ[q^{1/2},q^{-1/2}]$-basis
$\{T_w\theta_\lambda\mid w\in W, \lambda\in P\}$ of
$H(\tilde{W}_a)$ is due to Bernstein-Zelevinski (see \cite{LS}).

\section{affine Hecke algebras and equivariant $K$-groups}
\label{sec:HK} For an algebraic variety $Y$ over $\BC$ we denote
its structure sheaf by $\CO_Y$. If $Y$ is smooth, then its
canonical sheaf is denoted by $\Omega_Y$.

We denote the flag variety $G/B$ of $G$ by $\CB$. As a set $\CB$
is identified with the set of Borel subalgebras of $\Gg$ by the
correspondence $gB\mapsto\Ad(g)(\Gb)\,\,(g\in G)$. For $x\in\CB$
we denote by $\Gb_x$ the corresponding Borel subalgebra of $\Gg$,
and set $\Gn_x=[\Gb_x,\Gb_x]$. For $w\in W$ set
\[
Y_w= G(B,wB)\subset \CB\times\CB.
\]
Then we have $\CB\times\CB=\bigsqcup_{w\in W}Y_w$, and
$\overline{Y}_w=\bigsqcup_{y\leq w}Y_y$. We denote by
\[
i_w:\overline{Y}_w\to \CB\times\CB
\]
the embedding.

Set
\begin{align*}
\Lambda &=\{(a,x)\in\Gg\times\CB\mid a\in\Gn_x\},
\\
Z&= \{(a,x,y)\in\Gg\times\CB\times\CB\mid a\in\Gn_x\cap\Gn_y\}.
\end{align*}
Let $\pi:\Lambda\to\CB$ be the projection. The algebraic group
$G\times\BC^*$ acts on the variety $\Lambda$ by
\[
(g,z): \ (a,x)\mapsto(z^{-2}\Ad(g)(a),gx) \in\Lambda.
\]
We sometimes identify $Z$ with a $G\times\BC^*$-stable closed
subvariety of $\Lambda\times\Lambda$ by the embedding
\[
Z\to\Lambda\times\Lambda\qquad ((a,x,y)\mapsto((a,x),(a,y))).
\]
In particular, $Z$ is a $G\times\BC^*$-variety. For $w\in W$ set
\[
Z_w=\{(a,x,y)\in Z\mid (x,y)\in Y_w\}.
\]
We denote by
\[
r_w:\overline{Z}_w\to Z,\qquad \pi_w:\overline{Z}_w\to
\overline{Y}_w
\]
the embedding and the projection respectively.

Let us consider the equivariant $K$-group
$K^{G\times\BC^*}(Z)=K^{G\times\BC^*}(\Lambda\times\Lambda;Z)$
(see Section \ref{sec:K} for the equivariant $K$-groups and
notation concerning them). It is a module over the representation
ring $R^{G\times\BC^*}=R^G\otimes_\BZ R^{\BC^*}$ of $G\times\BC^*
$. We will identify $R^{\BC^*}$ with $\BZ[q^{1/2},q^{-1/2}]$ by
associating the $\BC^*$-module given by $z\mapsto z^n$ to
$q^{n/2}$. In particular, $K^{G\times\BC^*}(Z)$ is a
$\BZ[q^{1/2},q^{-1/2}]$-module.

For $(i,j)=(1,2), (2,3), (1,3)$ we denote by
$p_{ij}:\Lambda\times\Lambda\times\Lambda\to\Lambda\times\Lambda$
the projections onto $(i,j)$-factors. Note that
$p_{13}(p_{12}^{-1}Z\cap p_{23}^{-1}Z)\subset Z$. Since the
morphism $p_{12}^{-1}Z\cap p_{23}^{-1}Z\to Z$ induced by $p_{13}$
is proper, we can define an $R^{G\times\BC^*}$-bilinear map
\begin{align*}
& \bigstar :K^{G\times\BC^*}(Z)\times K^{G\times\BC^*}(Z)\to
K^{G\times\BC^*}(Z)
\\
&\quad \left( (m,n)\mapsto m\bigstar n=
p_{13*}(p_{12}^*m\otimes_{\CO_{\Lambda\times\Lambda\times\Lambda}}p_{23}
^*n) \right).
\end{align*}
Then it is easily seen that the convolution product $\bigstar$
endows with $K^{G\times\BC^*}(Z)$ a structure of  associative
algebra over $R^{G\times\BC^*}$ with the identity element
$[r_{1*}\CO_{Z_1}]$. For $\lambda\in P$ we denote by
$\CO_\CB(\lambda)$ the $G$-equivariant invertible $\CO_\CB$-module
whose fiber at $B$ is the $B$-module corresponding to $\lambda$.
\begin{theorem}
[Ginzburg \cite{G}, Kazhdan-Lusztig \cite{KL2}] \label{thm:HZ}
There exists an isomorphism
\[
\Phi:H(\tilde{W}_a)\to K^{G\times\BC^*}(Z)
\]
of  $\BZ[q^{1/2},q^{-1/2}]$-algebras satisfying
\begin{align*}
\Phi(\theta_\lambda)&= [r_{1*}\pi_1^*\CO_{\CB}(-\lambda)]
\qquad(\lambda\in P),\\
\Phi(T_s+1)&= -[r_{s*}\pi_s^*(\Omega_{\overline{Y}_s}\otimes
i_s^*(\CO_\CB\boxtimes\Omega_\CB^{\otimes-1}))] \qquad(s\in S).
\end{align*}
Here, we have identified $\overline{Y}_1(=Y_1)$ with $\CB$.
\end{theorem}
\begin{remark}
{\rm Note that $\Phi(T_s+1)$ is not symmetric with respect to the
the symmetry of $(\Lambda\times\Lambda,Z)$ given by
$\Lambda\times\Lambda\ni(x,y)\mapsto(y,x)\in\Lambda\times\Lambda$.
This can be resolved if we use the twisted product
\begin{align*}
&K^{G\times\BC^*}(Z)\times K^{G\times\BC^*}(Z)\to
K^{G\times\BC^*}(Z)
\\
&\quad \left((m,n)\mapsto
p_{13*}(p_{12}^*m\otimes_{\CO_{\Lambda\times\Lambda\times\Lambda}}p_{23}
^*n
\otimes_{\CO_{\Lambda\times\Lambda\times\Lambda}}p_2^*\pi^*\Omega_\CB)
\right)
\end{align*}
as in \cite{T}, where
$p_2:\Lambda\times\Lambda\times\Lambda\to\Lambda$ is the
projection onto the second factor. There is another way to recover
the symmetry by modifying the definition of $\Phi$ without
changing the product (see Lusztig \cite{L1}). }
\end{remark}
Let $\CN$ denote the closed subvariety of $\Gg$ consisting of
nilpotent elements. For a locally closed $G$-stable subvariety $V$
of $\CN$ we set
\[
Z_V=\{(a,x,y)\in Z\mid a\in V\}.
\]
\begin{proposition}
[Ginzburg \cite{G}, Kazhdan-Lusztig \cite{KL2}] \label{prop:exact}
Let $V$ be a locally closed $G$-stable subvariety of $\CN$. Then
we have an exact sequence
\[
0\to K^{G\times\BC^*}(Z_{\overline{V}\setminus V}) \to
K^{G\times\BC^*}(Z_{\overline{V}}) \to K^{G\times\BC^*}(Z_V) \to0.
\]
Here $K^{G\times\BC^*}(Z_{\overline{V}\setminus V}) \to
K^{G\times\BC^*}(Z_{\overline{V}})$ is given by the direct image
with respect to the inclusion $Z_{\overline{V}\setminus V} \to
Z_{\overline{V}}$, and $K^{G\times\BC^*}(Z_{\overline{V}}) \to
K^{G\times\BC^*}(Z_V)$ is given by the inverse image with respect
to the inclusion $Z_{V}\to Z_{\overline{V}}$.
\end{proposition}

In particular, if $V$ is closed, then the homomorphism
$K^{G\times\BC^*}(Z_V) \to K^{G\times\BC^*}(Z)$ given by the
direct image with respect to the closed embedding $Z_V\to Z$ is
injective. By this we will identify $K^{G\times\BC^*}(Z_V)$ for a
closed $G$-stable subvariety $V$ of $\CN$ with a two-sided ideal
of $K^{G\times\BC^*}(Z)$.

The following remarkable fact conjectured in Lusztig \cite{LC} was
proved  by Lusztig himself \cite{L:cell4} using the theory of
character sheaves among other things.
\begin{theorem}
\label{thm:bijection} There exists a natural one-to-one
correspondence between the set of two-sided cells of $\tilde{W}_a$
and that of nilpotent orbits of $\Gg$.
\end{theorem}
For a nilpotent orbit $O$ we denote by $\CC_O$ the corresponding
two-sided cell.

In view of Theorem \ref{thm:HZ}, it is natural to expect the
following (see \cite{G2, T2, L2}).
\begin{conjecture}
\label{conj:cell} Let $O$ be a nilpotent orbit. Then we have
\[
\Phi(H(\tilde{W}_a)_{\lrl  \CC_O})=
K^{G\times\BC^*}(Z_{\overline{O}}).
\]
\end{conjecture}
This conjecture is known to be true when $O=\{0\}$ (see
\cite{Xi}).

Let $w\in W$. In \cite{T} a $G\times\BC^{*}$-equivariant coherent
sheaf $M_w$ on $\Lambda\times\Lambda$ such that $\Supp(M_w)\subset
Z$ and $\Phi(C_w)=(-1)^{\ell(w)}[M_w]$ is associated using the
theory of Hodge modules. This together with a deep result related to the associated varieties of 
primitive ideals of the enveloping algebra $U(\Gg)$ implies
\[
\Phi(C_w) \in K^{G\times\BC^*}(Z_{\overline{O}}) \setminus
K^{G\times\BC^*}(Z_{\overline{O}\setminus O}),
\]
where $O$ is the nilpotent orbit satisfying $w\in\CC_O$. In
Section \ref{sec:GL} we will need the following weaker result
which is much easier.

\begin{proposition}
\label{prop:finite} Let $\Pi_1$ be a subset of $\Pi$. Set
$w=w_T\in W$ with $T=\{s_\alpha\mid\alpha\in\Pi_1\}\subset S$. Let
$O$ be the nilpotent orbit satisfying
\[
\overline{O}=
\Ad(G)(\sum_{\alpha\in\Delta^+\setminus\Delta_1}\Gg_\alpha),
\]
where $\Delta_1=\Delta\cap(\sum_{\alpha\in\Pi_1}\BZ\alpha)$. Then
we have
\[
\Phi(C_{w}) \in K^{G\times\BC^*}(Z_{\overline{O}}).
\]
\end{proposition}
\begin{proof}
Note that $\overline{Y}_w$ is smooth. Hence by \cite{T} we have
$\Phi(C_{w})=(-1)^{\ell({w})}[M_{w}]$ with
\[
M_w=\gr(\BQ^H_{\overline{Y}_w}[\dim
\overline{Y}_w])\otimes_{\CO_{\Lambda\times\Lambda}}(\pi\times\pi)^*(\CO
_\CB\boxtimes\Omega^{\otimes-1}_\CB),
\]
where $\BQ^H_{\overline{Y}_w}[\dim \overline{Y}_w]$ denotes the
canonical irreducible $G$-equivariant Hodge module whose
underlying perverse sheaf is $\BQ_{\overline{Y}_w}[\dim
\overline{Y}_w]$. By $\gr(\BQ^H_{\overline{Y}_w}[\dim
\overline{Y}_w]) =r_{w*}\pi_{w}^*(\Omega_{\overline{Y}_{w}})$ we
obtain
\begin{equation}
\label{eq:longest} M_{w}=
r_{w*}\pi_{w}^*(\Omega_{\overline{Y}_{w}})
\otimes_{\CO_{\Lambda\times\Lambda}}(\pi\times\pi)^*(\CO_\CB\boxtimes
\Omega^{\otimes-1}_\CB).
\end{equation}
It follows that $\Supp(M_{w})=\overline{Z}_{w}\subset
Z_{\overline{O}}$.
\end{proof}

\begin{remark}
{\rm We can prove \eqref{eq:longest} directly without appealing to
the theory of Hodge modules. Details are omitted. }
\end{remark}

\section{the case $G=GL_n(\BC)$}
\label{sec:GL} The main result of this paper is the following.
\begin{theorem}
\label{thm:main} Conjecture \ref{conj:cell} holds for
$G=GL_n(\BC)$.
\end{theorem}
In the rest of this section we assume that $G=GL_n(\BC)$. In this
case the extended affine Weyl group $\tilde{W}_a$ is identified
with the group of  all permutations $\sigma$ of $\BZ$ satisfying
$\sigma(i+n)=\sigma(i)+n\,\,(i\in\BZ)$ and
$\sum_{i=1}^n(\sigma(i)-i)\in n\BZ$. Define $\omega,
s_k\in\tilde{W}_a\,\,(0\leq k\leq n-1)$ by
\begin{align*}
{\omega}(i)&=i+1\qquad(i\in\BZ),\\
s_k(i)&=
\begin{cases}
i+1\qquad&(i\in n\BZ+k),\\
i-1\qquad&(i\in n\BZ+k+1),\\
i\qquad&(\text{otherwise}).
\end{cases}
\end{align*}
Then we have
\[
S=\{s_i\mid1\leq i\leq n-1\},\quad S_a=S\sqcup\{s_0\}, \quad
\Omega=\langle\omega\rangle,
\]
and $W$ is identified with the symmetric group $\mathfrak{S}_n$.

Let $\CP(n)$ denote the set of partitions of $n$, that is,
\[
\CP(n)= \{\rho=(\rho_1,\rho_2,\dots,\rho_n) \in\BZ_{\geq0}^n \mid
\rho_i\geq\rho_{i+1} ,\,\, \sum_{i=1}^n\rho_i=n \}.
\]
For $\rho\in\CP(n)$ we set
\[
N_j(\rho)= \sharp\{i\mid \rho_i=j\}.
\]
We denote by $\rho\mapsto\rho^*$ the duality operation on $\CP(n)$
induced by the transpose of the corresponding Young diagram, that
is, $\rho^*_i=\sum_{k=i}^nN_k(\rho)$.

The set of nilpotent orbits in $\Gg=\Ggl_n(\BC)$ is parametrized
by $\CP(n)$. The nilpotent orbit $O_\rho$ corresponding to
$\rho\in\CP(n)$ is the one containing the Jordan normal form with
exactly $N_i(\rho^*)$ Jordan blocks of  size $i$ (with eigenvalue
$0$) for each $i$. In particular, $O_{(n,0,\dots,0)}=\{0\}$ and
$O_{(1,\dots,1)}$ is the regular nilpotent orbit.

By Theorem \ref{thm:bijection} the set of two-sided cells of $\tilde W_a$ is also
parametrized by $\CP(n)$ (in our case $G=GL_n(\BC)$ this is due to
Lusztig \cite{L:cell} and Shi \cite{Shi}). We denote by $\CC_\rho$
the two-sided cell of $\tilde W_a$ corresponding to  $O_\rho$.

Let $T$ be a proper subset of $S_a$ such that $\langle T\rangle$
is of type $A_{k_1}\times\cdots\times A_{k_r}$. Then the
corresponding parabolic element $w_T$ belongs to $\CC_\rho$ if and
only if
\[
\sharp\{j\mid k_j+1=i\} =N_i(\rho)
\]
for any $i$.

For $\rho\in\CP(n)$ set $\CC_\rho^W=W\cap \CC_\rho$. It is known
that $\CC_\rho^W$ is a two-sided cell of $W$. In particular, the
set of two-sided cells of $W$ is also parametrized by $\CP(n)$
(see Kazhdan-Lusztig \cite{KL1}).
\begin{proposition}
[Shi \cite{Shi96}] \label{prop:key2} The following conditions on
$\rho, \xi\in\CP(n)$ are equivalent.
\begin{itemize}
\item[(a)] $\CC_\xi\lrl\CC_\rho$. \item[(b)]
$\CC^W_\xi\lrl\CC^W_\rho$. \item[(c)]
$O_\xi\subset\overline{O}_\rho$.
\end{itemize}
\end{proposition}
Hence we have $H(W)_{\lrl\CC_\rho^W}=H(\tilde
W_a)_{\underset{LR}\leq\CC_\rho}\cap H(W)$, and
$H(W)_{\CC^W_\rho}$ is identified with an $(H(W),H(W))$-submodule
of $H(\tilde W_a)_{\CC_\rho}$.

The following is crucial for the proof of Theorem \ref{thm:main}.

\begin{theorem}
\label{thm:key1} Let $v$ be a parabolic element in $\CC_\rho$.
Then the $H(\tilde{W}_a)$-bimodule $H(\tilde{W}_a)_{ \CC_\rho}$ is
generated by the image of $C_{v}$.
\end{theorem}

We first show the following corresponding statement for $H(W)$.
\begin{proposition}
\label{prop:key1} Let $v$ be a parabolic element in $\CC_\rho^W$.
Then the $H(W)$-bimodule $H(W)_{\CC_\rho^W}$ is generated by the
image of $C_{v}$.
\end{proposition}
\begin{proof}
Let $u\in\CC^W_\rho$. Let  $\CL$ be the left cell of $W$
containing $u$ and $\CR$ the right cell of $W$ containing $u$.
Then $\CL$ contains a unique element $y$ such that $y\er v$, and
$\CR$ contains a unique element $x$ such that $x\el v$
(see \cite{KL1}).
By Lemma \ref{lem:pa} we have $C'_x=hC'_v$ and $C'_y=C'_vh'$ for
some $h,h'$ in $H(W)$.

Let $\pi:H(W)_{\leq\CC^W_\rho}\to H(W)_{\CC^W_\rho}$ be the
canonical projection and let $V_1$ and $V_2$ be the left
$H(W)$-submodules of $H(W)_{\CC^W_\rho}$ generated by $\pi(C'_v)$
and $\pi(C'_y)$ respectively. Then $H(W)_{\CC^W_\rho}\ni k\mapsto
kh'\in H(W)_{\CC^W_\rho}$ is a homomorphism of left $H(W)$-modules
satisfying $\pi(C'_v)\mapsto \pi(C'_y)$. Hence we obtain a
homomorphism $f:V_1\to V_2$ given by $f(k)=kh'$.

On the other hand by \cite{KL1} there exists an isomorphism
$g:V_1\to V_2$ of left $H(W)$-modules such that
$g(\pi(C'_v))=\pi(C'_y)$ and $g(\pi(C'_x))=\pi(C'_u)$. By
$f(\pi(C'_v))=g(\pi(C'_v))$ we have $f=g$. Hence
\[
\pi(C'_u)=g(\pi(C'_x))=f(\pi(C'_x))=hf(\pi(C'_v)) =h\pi(C'_v)h'.
\]

The proof is complete.
\end{proof}

Now we give a proof  of Theorem \ref{thm:key1}. Let $\pi:H(\tilde
W_a)_{\leq\CC_\rho}\to H(\tilde W_a)_{\CC_\rho}$ be the canonical
homomorphism. According to \cite[Lemma 18.3.2]{Shi} one has a
parabolic element $w\in\CC^W_\rho$ such that for any
$u\in\CC_\rho$ there exists a sequence of left star operations
$\phi_1,\phi_2,...,\phi_r$ and an integer $m$ satisfying
\begin{equation}
\label{eq:star} w\er\omega^m\phi_r\phi_{r-1}\cdots\phi_1(u).
\end{equation}
We first show the statement for this special parabolic element
$w$.

Let $u\in\CC_\rho$. Take left star operations
$\phi_1,\phi_2,...,\phi_r$ and an integer $m$ satisfying
\eqref{eq:star}, and set
$y=\omega^m\phi_r\phi_{r-1}\cdots\phi_1(u)$,
$x=\phi_1\phi_2\cdots\phi_r\omega^{-m}(w)$. Note that $x$ is
well-defined and $x\el w$ by definition and Proposition
\ref{prop:LRstar}. Since $w$ is a parabolic element, there exist
$h, h'\in H(\tilde W_a)$ such that $C'_x=hC'_{w}$ and
$C'_y=C'_{w}h'$ by Lemma \ref{lem:pa}. Note that
$C'_{w}C'_{w}=\eta C'_{w}$ where $\eta\in\BZ[q^{1/2},q^{-1/2}]$
satisfies $\overline{\eta}=\eta$ and $\eta=q^{
{\ell(w)}/2}+\mbox{(lower degree terms)}$. Hence
\[
\eta
h\pi(C'_w)h'=\pi(C'_xC'_y)=\sum_{z\in\CC_\rho}h_{x,y,z}\pi(C'_z),
\]
where $h_{x,y,z}\in\BZ[q^{1/2},q^{-1/2}]$ satisfies
$\overline{h}_{x,y,z}=h_{x,y,z}$ and $h_{x,y,z}=\gamma_{x,y,z}q^{
{a(z)}/2}+\mbox{(lower degree terms)}$. For any $z\in\CC_\rho$ we
have $a(z)=a(w)=\ell(w)$, and hence we obtain
$$
h\pi(C'_{w})h'=\sum_{z\in\CC_\rho}\ga_{x,y,z}\pi(C'_z).$$ Note
that $\ga_{\omega^m w_1,w_2,w_3}=\gamma_{w_1,w_2,\omega^{-m}w_3}$
for any $w_1,w_2,w_3$ in $\tilde W_a$. Hence we have
$\gamma_{x,y,z}=\ga_{w,\ y,\ \omega^m\phi_r\cdots\phi_1(z)}$ by
Proposition \ref{prop:star}.  Since $w$ is a distinguished
involution, we have $\ga_{x,y,z}\ne 0$ if and only if
$\omega^m\phi_r\phi_{r-1}\cdots\phi_1(z)=y$ and in this case
$\ga_{x,y,z}=1$ (see Lusztig \cite{L:cell2}). Thus $\ga_{x,y,z}\ne
0$ if and only $z=u$ and in this case $\ga_{x,y,u}=1$. Therefore
we have $h\pi(C'_{w})h'=\pi(C'_u)$.

Now let $v$ be any parabolic element in $\CC_\rho$. Then there
exists an integer $k$ such that $\omega^kv\omega^{-k}$ is in $W$.
By Proposition \ref{prop:key1} we have
$H(W)\pi(C'_{\omega^kv\omega^{-k}})H(W)=H(W)\pi(C'_{w})H(W)$ and
hence
\[
H(\tilde W_a)\pi(C'_v)H(\tilde W_a) =H(\tilde
W_a)\pi(C'_{w})H(\tilde W_a) =H(\tilde W_a)_{\CC_\rho}.
\]
The proof of Theorem \ref{thm:key1} is complete.

\begin{remark}
{\rm (a) The assertion for $W_a$ similar to that in Theorem
\ref{thm:key1} does not hold in general.

(b) Let $v$ be as in Theorem \ref{thm:key1}. It is not difficult
to prove  that for any $w\llr v$, there exists a polynomial $f_w$
in $q^{\frac 12}+q^{-\frac12}$ such that $f_wC_w$ is in the
two-sided ideal of $H(\tilde{W}_a)$ generated by $C_{v}$. However,
in general it is  not true that $C_w$ is in the two-sided ideal of
$H(\tilde{W}_a)$ generated by $C_{v}$. Example: $n=4$ and let
$\CC_\rho$ be the two-sided cell containing $v=s_1s_3$. Then
$(q^{\frac 12}+q^{-\frac12})C_{s_1s_2s_1}$ is in
$H(\tilde{W}_a)C_{v}H(\tilde{W}_a)$, but $C_{s_1s_2s_1}$ is not in
$H(\tilde{W}_a)C_{v}H(\tilde{W}_a)$. }
\end{remark}

Let $\BF$ be an algebraic closure of $\BC(q^{1/2})$, and set
$H^\BF=\BF\otimes H(\tilde{W}_a), G_\BF=GL_n(\BF),
\Gg_\BF=\Ggl_n(\BF)$. Then $H^\BF$ is an $\BF$-algebra and $G_\BF$
is an algebraic group over $\BF$ with Lie algebra $\Gg_\BF$.

Let $\CQ$ denote the $G_\BF$-conjugacy classes of the pairs
$(s,e)\in G_{\BF}\times\Gg_\BF$ where $s$ is semisimple,  $e$ is
nilpotent, and $\Ad(s)(e)=qe$. For  such a pair $(s,e)$
Kazhdan-Lusztig \cite{KL2} and Ginzburg \cite{G} constructed a
finite-dimensional $H^\BF$-module $M_{(s,e)}$. Moreover, we have a
unique irreducible quotient $L_{(s,e)}$ of $M_{(s,e)}$, and the
set of irreducible $H^\BF$-modules is parametrized by  $\CQ$ via
$(s,e)\mapsto L_{(s,e)}$ (note that $\BF$ is isomorphic to $\BC$
as an abstract field). In particular, we can associate to each
irreducible $H^\BF$-module $L$ a nilpotent orbit $O(L)$ in $\Gg$
by $\Ad(G_\BF)(O(L_{(s,e)}))\ni e$ (note that the set of
$G_\BF$-conjugacy classes of nilpotent elements in $\Gg_\BF$ is in
one-to-one correspondence with that of $G$-conjugacy classes of
nilpotent elements in $\Gg$).

We need the following deep result of Lusztig \cite{L:cell4}.
\begin{proposition}
\label{prop:lusztig-cell} For any irreducible subquotient $L$ of
the $($left$)$ $H^\BF$-module
$\BF\otimes_{\BZ[q^{1/2},q^{-1/2}]}H(\tilde{W}_a)_{\CC_\rho}$ we
have $\overline{O(L)}\supset O_\rho$.
\end{proposition}
\begin{proposition}
\label{prop:FK} Let $O$ be a nilpotent orbit. Then for any
irreducible quotient $L$ of the $($left$)$ $H^\BF$-module
$\BF\otimes_{\BZ[q^{1/2},q^{-1/2}]}K^{G\times\BC^*}(Z_{O})$ we
have $O(L)=O$.
\end{proposition}
\begin{proof}
By \cite[Corollary 5.9]{KL2} we see that $L$ is a quotient of
$M_{(s,e)}$ for $(s,e)\in\CQ$ with $e\in O$. Since $L_{(s,e)}$ is
the unique irreducible quotient of $M_{(s,e)}$, we have
$L=L_{(s,e)}$ and hence $O(L)=O$.
\end{proof}

Now we are ready to give a proof of Theorem \ref{thm:main}. We
show
\begin{equation}
\label{eq:claim} \Phi(H(\tilde{W}_a)_{\llr\CC_\xi})=
K^{G\times\BC^*}(Z_{\overline{O}_\xi})
\end{equation}
for any $\xi\in\CP(n)$ by induction on $\dim O_\xi$. Let
$\rho\in\CP(n)$ and assume that \eqref{eq:claim} is true for any
$\xi\in\CP(n)$ with $\dim O_\xi<\dim O_\rho$.

For any $\tau\in\CP(n)$ with $\CC_\tau\llr\CC_\rho$ any parabolic
element $v\in\CC_\tau^W$ satisfies $\Phi(C_{v})\in
K^{G\times\BC^*}(Z_{\overline{O}_\rho})$ by Proposition
\ref{prop:finite}. Hence we see by Theorem \ref{thm:key1} that
$\Phi(H(\tilde{W}_a)_{\llr \CC_\rho}) \subset
K^{G\times\BC^*}(Z_{\overline{O}_\rho})$. Moreover, the hypothesis
of induction together with Proposition \ref{prop:key2} implies
$\Phi(H(\tilde{W}_a)_{\underset{LR}< \CC_\rho}) = K^{G\times\BC^*}
(Z_{\overline{O}_\rho\setminus {O}_\rho})$. Hence it is sufficient
to show that the induced injection
$\overline{\Phi}:H(\tilde{W}_a)_{\CC_\rho} \to
K^{G\times\BC^*}(Z_{{O}_\rho})$ is surjective. Assume that
$\Coker(\overline{\Phi})\ne0$. Since  $H(\tilde{W}_a)_{\llr
\CC_\rho}$  is a direct summand of the
$\BZ[q^{1/2},q^{-1/2}]$-module $H(\tilde{W}_a)$ and $\Phi(H(\tilde
W_a))=K^{G\times C^*}(Z)$, we see that the cokernel of  the
injective homomorphism
\[
\overline{\Phi}^\BF: \BF\otimes_{\BZ[q^{1/2},q^{-1/2}]}
H(\tilde{W}_a)_{\CC_\rho} \to \BF\otimes_{\BZ[q^{1/2},q^{-1/2}]}
K^{G\times\BC^*}(Z_{{O}_\rho})
\]
is also non-trivial. Take an irreducible quotient $L$ of the
$H^\BF$-module $\Coker(\overline{\Phi}^\BF)$. Since $L$ is an
irreducible quotient of $\BF\otimes_{\BZ[q^{1/2},q^{-1/2}]}
K^{G\times\BC^*}(Z_{{O}_\rho})$, we have $O(L)=O_\rho$ by
Proposition  \ref{prop:FK}. On the other hand since $L$ is an
irreducible subquotient of the $H^\BF$-module
$$\BF\otimes_{\BZ[q^{1/2},q^{-1/2}]} H(\tilde{W}_a)/
\BF\otimes_{\BZ[q^{1/2},q^{-1/2}]}
H(\tilde{W}_a)_{\llr\CC_\rho},$$ there exists a nilpotent orbit
$O$ such that $O\not\subset\overline{O_\rho}$ and
$O\subset\overline{O(L)}$ by Proposition \ref{prop:lusztig-cell}.
This is a contradiction. Hence $\overline{\Phi}$ is surjective.
The proof of  Theorem \ref{thm:main} is complete.

\appendix

\section{Equivariant $K$-theory}
\label{sec:K} In this section we recall basic notions concerning
equivariant $K$-groups (see Thomason \cite{Th}). All algebraic
varieties are assumed to be quasi-projective over $\BC$ and all
algebraic groups are assumed to be affine over $\BC$. The
structure sheaf of an algebraic variety $X$ is denoted by $\CO_X$.
When we consider an action of an algebraic group $A$ on an
algebraic variety $X$, we always assume the existence of a closed
$A$-equivariant embedding $X\to X'$ where $X'$ is a smooth variety
with an action of $A$. In this case we say that $X$ is an
$A$-variety.

Let $A$ be an algebraic group. For a pair $(Y,X)$ such that $Y$ is
an $A$-variety and $X$ is its $A$-stable closed subvariety, we
denote by $\Coh^A(Y;X)$ the abelian category of  $A$-equivariant
coherent sheaves on $Y$ whose supports are contained in $X$. Its
Grothendieck group $K^A(Y;X)$ is called the equivariant $K$-group.
Note that the direct image functor $i_*:\Coh^A(X;X)\to\Coh^A(Y;X)$
with respect to the embedding $i:X\to Y$ induces an isomorphism
$K^A(X;X)\cong K^A(Y;X)$. It means that $K^A(Y;X)$ depends only on
the $A$-variety $X$, and hence we sometimes denote it by $K^A(X)$.
However, we will need to specify the ambient space $Y$ in defining
some operations on equivariant $K$-groups. Note that $K^A(X)$ is a
module over the representation ring
\begin{equation}
\label{eq:K:1} R^A=K^A(\poi)
\end{equation}
of $A$. Here $\poi$ denotes the variety consisting of a single
point.

Assume that we are given an $A$-equivariant morphism $f:Y\to Y'$
of $A$-varieties and $A$-stable closed subvarieties $X$ and $X'$
of $Y$ and $Y'$ respectively such that $f(X)\subset X'$ and the
restriction $X\to X'$ of $f$ is a proper morphism. Then the
derived functors
\[
R^nf_*:\Coh^A(Y;X)\to\Coh^A(Y';X')\qquad(n\in\BZ)
\]
of the direct image functor $f_*$ induce a homomorphism
\begin{equation}
\label{eq:K:2} f_*:K^A(Y;X)\to K^A(Y';X')\quad
([M]\mapsto\sum_n(-1)^n[R^nf_*(M)])
\end{equation}
of $R^A$-modules. We note that \eqref{eq:K:2} does not depend on
the choice of the ambient spaces $Y$ and $Y'$.
\begin{lemma}
Let $f:X\to X'$ and $g:X'\to X''$ be $A$-equivariant proper
morphisms of $A$-varieties. Then we have
\[
(g\circ f)_*=g_*\circ f_*:K^A(X)\to K^A(X'').
\]
\end{lemma}

Assume that we are given an $A$-equivariant morphism $f:Y\to Y'$
of $A$-varieties and an $A$-stable closed subvariety $X'$ of $Y'$.
Set $X=f^{-1}(X')$. If $f$ is smooth or if $Y'$ is a smooth
variety, then the derived functors
\[
L^nf^*:\Coh^A(Y';X')\to\Coh^A(Y;X)\qquad(n\in\BZ)
\]
of  the inverse image functor
\[
f^*:\Coh^A(Y';X')\to\Coh^A(Y;X)\quad (M\mapsto
\CO_Y\otimes_{f^{-1}\CO_{Y'}}f^{-1}M)
\] are zero except for finitely many $n$'s, and they induce a
homomorphism
\begin{equation}
\label{eq:K:3} f^*:K^A(Y';X')\to K^A(Y;X)\quad
([M]\mapsto\sum_n(-1)^n[L^nf^*(M)])
\end{equation}
of $R^A$-modules.
  If $f$ is smooth, we have $L^nf^*=0$ for $n\ne0$.
  \begin{lemma}
Let $f:Y\to Y'$ and $g:Y'\to Y''$ be $A$-equivariant morphisms of
$A$-varieties. Let $X''$ be a closed subvariety of $Y''$, and set
$X=(g\circ f)^{-1}(X''), X'=f^{-1}(X'')$. Assume that
$f^*:K^A(Y';X')\to K^A(Y;X)$ and $g^*:K^A(Y'';X'')\to K^A(Y';X')$
are defined. Then we have
\[
(g\circ f)^*=f^*\circ g^*: K^A(Y'';X'')\to K^A(Y;X).
\]
\end{lemma}

Assume that we are given a smooth $A$-variety  $Y$ and its
$A$-stable closed subvarieties $X_1$ and $X_2$. The derived
functors
\begin{align*}
\Tor_n^{\CO_Y}(\,\,,\,\,): \Coh^A(Y;X_1)\times\Coh^A(Y;X_2)\to
\Coh^A(Y;X_1\cap X_2)\\
((M_1,M_2)\mapsto\Tor_n^{\CO_Y}(M_1,M_2)=
H^{-n}(M_1\otimes_{\CO_Y}^{\BLL}M_2))
\end{align*}
of the tensor product functor $\otimes_{\CO_Y}$ are zero except
for finitely many $n$'s, and induce a bilinear map
\begin{align}
\label{eq:K:4} \otimes_{\CO_Y}: K^A(Y;X_1)\times K^A(Y;X_2)\to
K^A(Y;X_1\cap X_2)\\
(([M_1],[M_2])\mapsto [M_1]\otimes_{\CO_Y}[M_2]=
\sum_n(-1)^n\Tor_n^{\CO_Y}(M_1,M_2)) \nonumber
\end{align}
of  $R^A$-modules. Note that $\otimes_{\CO_Y}$ does depend on the
choice of the ambient space $Y$.
\begin{lemma}\label{lem:tens}
Let $f:Y\to Y'$ be an $A$-equivariant smooth morphism of smooth
$A$-varieties. Let $X'_1, X'_2$ be  closed subvarieties of $Y'$,
and set $X_1=f^{-1}(X'_1), X_2=f^{-1}(X'_2)$. Then we have
\[
f^*(m_1)\otimes_{\CO_Y}f^*(m_2) =f^*(m_1\otimes_{\CO_{Y'}}m_2) \in
K^A(Y;X_1\cap X_2)
\]
for any $m_1\in K^A(Y';X_1'), m_2\in K^A(Y';X_2')$.
\end{lemma}
\begin{lemma}[projection formula]
\label{lem:pr} Let $f:Y\to Y'$ be an $A$-equivariant morphism of
smooth $A$-varieties. Let $X_1'$ be an $A$-stable closed
subvariety of $Y'$ and set $X_1=f^{-1}(X'_1)$. Let $X_2$ and
$X'_2$ be  closed subvarieties of $Y$ and $Y'$ respectively such
that $f(X_2)=X_2'$ and $X_2\to X'_2$ is proper. Then we have
\[
f_*(f^*(m)\otimes_{\CO_Y}n) =m\otimes_{\CO_{Y'}}f_*n \in
K^A(Y';X_1'\cap X_2')
\]
for any $m\in K^A(Y';X_1'), n\in K^A(Y;X_2)$.
\end{lemma}
\begin{lemma}
[base change theorem 1] \label{lem:BC1} Let $f:Y'\to Y$ and
$g:Y''\to Y$ be $A$-equivariant morphism of $A$-varieties. We
assume that $g$ is smooth. Set $Y'''=Y'\times_YY''$ and let
$f':Y'''\to Y''$ and $g':Y'''\to Y'$ be canonical morphisms. Let
$X, X'$ be closed $A$-stable closed subvarieties of $Y, Y'$
respectively such that $f(X')\subset X$ and $X'\to X$ is proper.
Then we have
\[
g^*\circ f_*=f'_*\circ g'^*:K^A(Y';X')\to K^A(Y'';g^{-1}(X)).
\]
\end{lemma}
\begin{lemma}[base change theorem 2]
\label{lem:BC2} Let $Y$ be a smooth $A$-variety and let $Y_1, Y_2$
be $A$-stable smooth closed subvarieties of $Y$. Set $Y_3=Y_1\cap
Y_2$. We assume that $Y_3$ is smooth and that
\[
T_yY=T_yY_1+T_yY_2,\qquad T_yY_3=T_yY_1\cap T_yY_2
\]
for any $y\in Y_3$. Here, $T_yY$ denotes the tangent space of $Y$
at $y$. Let $i:Y_1\to Y, \,\, j:Y_2\to Y,\,\,  i':Y_3\to Y_2,\,\,
j':Y_3\to Y_1$ be the inclusions. Let $X_1$ be an $A$-stable
closed subvariety of $Y_1$. Then we have
\[
j^*\circ i_*=i'_*\circ j'^*:K^A(Y_1;X_1) \to K^A(Y_2;X_1\cap Y_2).
\]
\end{lemma}
\section{Convolution product}
\label{sec:conv} In this section $G$ is as in  Section
\ref{sec:affine Hecke}. In particular, $G$ is not necessarily of
type $A$. We fix a nilpotent orbit $O$ of $\Gg$ in the following.

According to Conjecture \ref{conj:cell} the quotient
\[
H(\tilde{W}_a)_{\CC_O}= H(\tilde{W}_a)_{\lrl
\CC_O}/H(\tilde{W}_a)_{\underset{LR}<\CC_O}
\]
should be identified with
\[
K^{G\times\BC^*}(Z_{{O}})\cong K^{G\times\BC^*}(Z_{\overline{O}})/
K^{G\times\BC^*}(Z_{\overline{O}\setminus O}).
\]
For $e\in O$ set
\[
\CB_e=\{x\in\CB\mid e\in\Gn_x\}.
\]
Since $Z_O$ is a $G\times\BC^*$-equivariant fiber bundle on $O$
whose fiber at $e\in O$ is canonically isomorphic to
$\CB_e\times\CB_e$, we have
\begin{equation}
\label{eq:conv:1} K^{G\times\BC^*}(Z_{{O}}) \cong
K^{M(e)}(\CB_e\times\CB_e),
\end{equation}
where
\[
M(e) =\{(g,z)\in G\times\BC^*\mid \Ad(g)(e)=z^2e\}.
\]

The aim of this section is to give a description of the product on
$K^{M(e)}(\CB_e\times\CB_e)$ induced from the convolution product
$\bigstar$ on $K^{G\times\BC^*}(Z)$.

We say that a triple $(h, e, f)\in\Gg\times\Gg\times\Gg$ is an
$\Gsl_2$-triple if $[h,e]=2e, [h,f]=-2f, [e,f]=h$. Then $e$ and
$f$ are nilpotent elements belonging to the same conjugacy class.
Moreover, the map $(h,e,f)\mapsto e$ induces a bijection between
the set of $G$-conjugacy classes of  $\Gsl_2$-triples and that of
nilpotent orbits. Set
\[
\hat{O} =\{(e,f)\in\Gg\times\Gg\mid e\in O, \text{ $([e,f],e,f)$
is an $\Gsl_2$-triple} \}.
\]
The group $G\times\BC^*$ acts transitively on $\hat{O}$ by
\[
(g,z):(e,f)\to(z^{-2}\Ad(g)(e),z^{2}\Ad(g)(f)).
\]
In particular, $\hat{O}$ is a smooth variety. For
$(e,f)\in\hat{O}$, Slodowy's variety $\Lambda_{(e,f)}$ is defined
by
\[
\Lambda_{(e,f)}= \{(a,x)\in\Lambda\mid a\in e+\Gz_\Gg(f)\},
\]
where
\[
\Gz_\Gg(f)= \{a\in\Gg\mid[a,f]=0\}.
\]
\begin{proposition}
[Slodowy \cite{S}] \label{prop:Slodowy}
\begin{itemize}
\item[(i)] $\Lambda_{(e,f)}$ is a smooth variety with
$\dim\Lambda_{(e,f)}=2\dim\CB_e$. \item[(ii)]
$\Ad(G)(\CN\cap(e+\Gz_\Gg(f))\subset\CN\setminus(\overline{O}\setminus
O)$. \item[(iii)] $O\cap(e+\Gz_\Gg(f))=\{e\}$.
\end{itemize}
\end{proposition}
We identify $\CB_e$ with a closed subvariety of $\Lambda_{(e,f)}$
via the embedding $x\mapsto(e,x)$. Set
\[
M(e,f) =\{(g,z)\in G\times\BC^*\mid \Ad(g)(e)=z^{2}e,\,\,
\Ad(g)(f)=z^{-2}f\}.
\]
Then $M(e,f)$ is a subgroup of $M(e)$ acting naturally on
$\Lambda_{(e,f)}$. Moreover, $M(e,f)$ and $M(e)$ contain a common
maximal reductive subgroup (see \cite{KL2}). Hence we have the
identification
\begin{equation}
\label{eq:conv:2} K^{M(e)}(\CB_e\times\CB_e) =
K^{M(e,f)}(\Lambda_{(e,f)}\times\Lambda_{(e,f)};\CB_e\times\CB_e).
\end{equation}
For $(i,j)=(1,2), (2,3), (1,3)$ we denote by $\pi_{ij}:
\Lambda_{(e,f)}\times\Lambda_{(e,f)}\times\Lambda_{(e,f)}
\to\Lambda_{(e,f)}\times\Lambda_{(e,f)}$ the projections onto
$(i,j)$-factors.
\begin{theorem}
\label{thm:conv} The product on
$K^{M(e,f)}(\Lambda_{(e,f)}\times\Lambda_{(e,f)};\CB_e\times\CB_e)$
induced from the convolution product $\bigstar$ on
$K^{G\times\BC^*}(Z)$ is given by
\[
(m,n)\mapsto
\pi_{13*}(\pi_{12}^*m\otimes_{\CO_{\Lambda_{(e,f)}\times\Lambda_{(e,f)}
\times\Lambda_{(e,f)}}}\pi_{23}^*n).
\]
\end{theorem}
The rest of this section is devoted to proving Theorem
\ref{thm:conv}.

Set
\[
\tilde{\Lambda} =\{(a,x)\in\Lambda\mid
a\notin\overline{O}\setminus O\}.
\]
Then $\tilde{\Lambda}$ is an open subset of $\Lambda$, and $Z_O$
is a closed subset of $\tilde{\Lambda}\times\tilde{\Lambda}$. We
denote by
\[
k:Z_O\to\tilde{\Lambda}\times\tilde{\Lambda}
\]
the closed embedding. For $(i,j)=(1,2), (2,3), (1,3)$ we denote by
$\tilde{p}_{ij}:
\tilde{\Lambda}\times\tilde{\Lambda}\times\tilde{\Lambda}\to\tilde
{\Lambda} \times\tilde{\Lambda}$ the projections onto
$(i,j)$-factors. We see easily the following.
\begin{lemma}
\label{lem:open} The product on $K^{G\times\BC^*}(Z_{{O}})
=K^{G\times\BC^*}(\tilde{\Lambda}\times\tilde{\Lambda};Z_{{O}})$
induced from the convolution product $\bigstar$ on
$K^{G\times\BC^*}(Z)$ is given by
\[
(m,n)\mapsto m\bigstar n=
\tilde{p}_{13*}(\tilde{p}_{12}^*m\otimes_{\CO_{\tilde{\Lambda}\times\tilde
{\Lambda}\times\tilde{\Lambda}}}\tilde{p}_{23}^*n).
\]
\end{lemma}
Set
\begin{align*}
\tilde{\Lambda}_O &= \{(e,x)\in\Lambda\mid
e\in O\},\\
Y_O &= \tilde{\Lambda}_O\times_O\hat{O} =\{(e,f,x)\mid
(e,f)\in\hat{O},\,\,(e,x)\in\Lambda\},\\
Y &= \{(e,f,a,x)\mid
(e,f)\in\hat{O},\,\,(a,x)\in\Lambda_{(e,f)}\}.
\end{align*}
We identify $Y_O$ with a closed subvariety of $Y$ by the embedding
\[
i:Y_O\to Y\qquad ((e,f,x)\mapsto(e,f,e,x)).
\]
Then $Y$ is a $G\times\BC^*$-equivariant fiber bundle on $\hat{O}$
whose fiber at $(e,f)\in\hat{O}$ is $\Lambda_{(e,f)}$, and $Y_O$
is its subbundle whose fiber at $(e,f)\in\hat{O}$ is $\CB_e$. In
particular, $Y$ is a smooth variety and the projection
$Y\to\hat{O}$ is a smooth morphism.

We set
\begin{align*}
Y^{(2)}&=Y\times_{\hat{O}}Y\\
&=\{(e,f,a,x,b,y)\mid (e,f)\in\hat{O},\,\,
(a,x), (b,y)\in\Lambda_{(e,f)}\},\\
Y_O^{(2)}&=Y_O\times_{\hat{O}}Y_O
=Z_O\times_O\hat{O}\\
&=\{(e,f,x,y)\mid (e,f)\in\hat{O},\,\, x,y\in\CB_e\}.
\end{align*}
We regard $Y_O^{(2)}$ as a closed subvariety of $Y^{(2)}$ by the
embedding
\[
i^{(2)}=i\times_{\hat{O}}i:Y_O^{(2)}\to Y^{(2)}.
\]
Define
\[
\varphi:Y_O^{(2)}\to Z_O
\]
by $\varphi(e,f,x,y)=(e,x,y)$. It is a smooth surjective morphism.
Since $Y_O^{(2)}$ is a $G\times\BC^*$-equivariant fiber bundle
whose fiber at $(e,f)\in\hat{O}$ is
$\Lambda_{(e,f)}\times\Lambda_{(e,f)}$, we have a commutative
diagram
\[
\begin{CD}
K^{G\times\BC^*}(Z_O) @>{\varphi^*}>>
K^{G\times\BC^*}(Y_O^{(2)})\\
@|@|\\
K^{M(e)}(\CB_e\times\CB_e) @>>> K^{M(e,f)}(\CB_e\times\CB_e).
\end{CD}
\]
Hence we see by \eqref{eq:conv:2} that
\begin{equation}
\label{eq:conv:3} \varphi^*:K^{G\times\BC^*}(Z_O) \to
K^{G\times\BC^*}(Y_O^{(2)})
\end{equation}
is an isomorphism of $R^{G\times\BC^*}$-modules. Set
\begin{align*}
Y^{(3)}&=Y\times_{\hat{O}}Y\times_{\hat{O}}Y,\\
Y_O^{(3)}&=Y_O\times_{\hat{O}}Y\times_{\hat{O}}Y,
\end{align*}
and regard $Y_O^{(3)}$ as a subvariety of $Y^{(3)}$ by
\[
i^{(3)}=i\times_{\hat{O}}i\times_{\hat{O}}i :Y_O^{(3)}\to Y^{(3)}.
\]
For $(i,j)=(1,2), (2,3), (1,3)$ we denote by $q_{ij}:Y^{(3)}\to
Y^{(2)}$ the projections onto $(i,j)$-factors. Note that $q_{ij}$
is a morphism of $G\times\BC^*$-equivariant fiber bundles on
$\hat{O}$ whose fiber at $(e,f)\in\hat{O}$ is given by
$\pi_{ij}:\Lambda_{(e,f)}\times\Lambda_{(e,f)}\times\Lambda_{(e,f)}
\to\Lambda_{(e,f)}\times\Lambda_{(e,f)}$. Therefore, Theorem
\ref{thm:conv} is equivalent to the following.
\begin{proposition}
\label{prop:conv} The product on
$K^{G\times\BC^*}(Y^{(2)};Y_O^{(2)})$ induced from the convolution
product $\bigstar$ on $K^{G\times\BC^*}(Z_O)$ via $\varphi^*$ is
given by
\[
(m,n)\mapsto q_{13*}(q_{12}^*m\otimes_{\CO_{Y^{(3)}}}q_{23}^*n).
\]
\end{proposition}
By Proposition \ref{prop:Slodowy} we have a morphism
\[
\theta:Y\to\tilde{\Lambda} \qquad((e,f,a,x)\mapsto(a,x)).
\]
We define
\[
\tau:Y_O\to\tilde{\Lambda}_O
\]
as the restriction of $\theta$.
\begin{lemma}
\label{lem:conv:theta}
\begin{itemize}
\item[(i)] The commutative diagram
\[
\begin{CD}
Y_O@>{\tau}>>\tilde{\Lambda}_O\\
@VVV@VVV\\
Y@>>{\theta}>\tilde{\Lambda}
\end{CD}
\]
is cartesian. \item[(ii)] $\theta$ is a smooth morphism.
\end{itemize}
\end{lemma}
\begin{proof}
The statement (i) follows from Proposition \ref{prop:Slodowy}
(iii). By a result of Slodowy \cite{S} we see that the composition
of the smooth surjective morphism
\[
G\times \Lambda_{(e,f)}\to Y\qquad((g,(a,x))\mapsto
(\Ad(g)(e),\Ad(g)(f),\Ad(g)(a),gx))
\]
with $\theta:Y\to \tilde{\Lambda}$ is smooth (see the proof of
Proposition 11.10 in Lusztig \cite{L1}). Hence $\theta$ is also
smooth.
\end{proof}

Consider the following diagrams
\begin{align}
\label{diagram1} &
\begin{CD}
Y_O^{(2)}\times_{\hat{O}}Y @>{i^{(2)}\times_{\hat{O}}1}>> Y^{(3)}
@>{k_{23}}>> \tilde{\Lambda}\times Y^{(2)}
\\
@V\beta_{12}VV @VV{k_{12}}V @VV{\ell_{23}}V
\\
Y_O^{(2)}\times\tilde{\Lambda} @>>{i^{(2)}\times1}>
Y^{(2)}\times\tilde{\Lambda} @>>\ell_{12}> \tilde{\Lambda}\times
Y\times\tilde{\Lambda}
\\
@V{\alpha_{12}}VV @V{\gamma_{12}}VV
\\
Y_O^{(2)} @>>{i^{(2)}}> Y^{(2)},
\\
\end{CD}
\\
\nonumber
\\
\label{diagram2} &
\begin{CD}
Y_O^{(2)}\times\tilde{\Lambda}
@>{\ell_{12}\circ(i^{(2)}\times1)}>> \tilde{\Lambda}\times
Y\times\tilde{\Lambda}
\\
@V{\varphi\circ\alpha_{12}}VV
@VV{\tilde{p}_{12}\circ(1\times\theta\times1)}V
\\
Z_O @>>{k}> \tilde{\Lambda}\times\tilde{\Lambda}
\end{CD}
\end{align}
where $\alpha_{12}, \gamma_{12}$ are the projections, and
$\beta_{12}, k_{12}, k_{23}, \ell_{12}, \ell_{23}$ are the closed
embeddings  induced by $\theta:Y\to \tilde{\Lambda}$. We can check
the commutativity easily. Moreover, we see easily that all of the
squares in the diagrams are cartesian. We set
\[
\psi=\ell_{12}\circ k_{12} =\ell_{23}\circ k_{23}:Y^{(3)}\to
\tilde{\Lambda}\times Y\times\tilde{\Lambda}.
\]

Let $m, n\in K^{G\times\BC^*}(Z_O;Z_O)$. Then the corresponding
elements in $K^{G\times\BC^*}(Y^{(2)};Y_O^{(2)})$ are given by
$\tilde{m}=i^{(2)}_*\varphi^*m,\,\, \tilde{n}=i^{(2)}_*\varphi^*n$
respectively. By $q_{12}=\gamma_{12}\circ k_{12}$ we see  from
\eqref{diagram1} that
\[
q_{12}^*\tilde{m} =k_{12}^*\gamma_{12}^*i^{(2)}_*\varphi^*m
=k_{12}^*(i^{(2)}\times1)_*\alpha_{12}^* \varphi^*m.
\]
Similarly, we have $q_{23}^*\tilde{n}= k_{23}^*(1\times
i^{(2)})_*\alpha_{23}^* \varphi^*n$, where
$\alpha_{23}:\tilde{\Lambda}\times Y_O^{(2)}\to Y_O^{(2)}$ is the
projection. Hence we have
\begin{align*}
&\psi_*(q_{12}^*\tilde{m} \otimes
q_{23}^*\tilde{n})\\
=&\ell_{12*}k_{12*}(k_{12}^*(i^{(2)}\times1)_*\alpha_{12}^*
\varphi^*m \otimes k_{23}^*(1\times i^{(2)})_*\alpha_{23}^*
\varphi^*n)\\
=& \ell_{12*}((i^{(2)}\times1)_*\alpha_{12}^* \varphi^*m \otimes
k_{12*}k_{23}^*(1\times i^{(2)})_*\alpha_{23}^*
\varphi^*n)\\
=& \ell_{12*}((i^{(2)}\times1)_*\alpha_{12}^* \varphi^*m \otimes
\ell_{12}^*\ell_{23*}(1\times i^{(2)})_*\alpha_{23}^*
\varphi^*n)\\
=& \ell_{12*}(i^{(2)}\times1)_*\alpha_{12}^* \varphi^*m \otimes
\ell_{23*}(1\times i^{(2)})_*\alpha_{23}^* \varphi^*n.
\end{align*}
Here we have used Lemma \ref{lem:pr} for the second and the fourth
identities and Lemma  \ref{lem:BC2} for the third identity. By
Lemma  \ref{lem:BC1}, Lemma \ref{lem:conv:theta} and
\eqref{diagram2} we have
\[
\ell_{12*}(i^{(2)}\times1)_*\alpha_{12}^* \varphi^*m
=(1\times\theta\times1)^*\tilde{p}_{12}^*k_*m.
\]
Similarly we have
\[
\ell_{23_*}(1\times i^{(2)})_*\alpha_{23}^* \varphi^*n
=(1\times\theta\times1)^*\tilde{p}_{23}^*k_*n.
\]
Therefore, we obtain
\[
\psi_*(q_{12}^*\tilde{m} \otimes q_{23}^*\tilde{n}) =
(1\times\theta\times1)^* (\tilde{p}_{12}^*k_*m \otimes
\tilde{p}_{23}^*k_*n)
\]
by Lemma \ref{lem:tens}.

Set
\[
\tilde{\Lambda}_O^{(3)}= \tilde{\Lambda}_O\times_O
\tilde{\Lambda}_O\times_O \tilde{\Lambda}_O
=\tilde{p}_{12}^{-1}Z_O\cap\tilde{p}_{23}^{-1}Z_O,
\]
and consider the commutative diagram
\begin{equation}
\label{diagram3}
\begin{CD}
\tilde{\Lambda}\times\tilde{\Lambda}\times\tilde{\Lambda}
@<{1\times\theta\times 1}<< \tilde{\Lambda}\times
Y\times\tilde{\Lambda}
\\
@A{f}AA @AA{{\psi}\circ i^{(3)}}A
\\
\tilde{\Lambda}_O^{(3)} @<{\tilde{\varphi}}<< Y_O^{(3)}
\\
@V{\overline{p}_{13}}VV @VV{\overline{q}_{13}}V
\\
Z_O @<<{\varphi}< Y_O^{(2)}.
\end{CD}
\end{equation}
Here, $f$ is the natural inclusion, $\tilde{\varphi}$ is the
canonical morphism, and $\overline{p}_{13}, \overline{q}_{13}$ are
the restrictions of $\tilde{p}_{13}, {q}_{13}$ respectively. We
see easily that both of the squares in \eqref{diagram3} are
cartesian.

Define $u\in K^{G\times\BC^*} (\tilde{\Lambda}_O^{(3)}
;\tilde{\Lambda}_O^{(3)})$ by $f_*u= \tilde{p}_{12}^*k_*m \otimes
\tilde{p}_{23}^*k_*n$. Then we have
\[
\psi_*(q_{12}^*\tilde{m} \otimes q_{23}^*\tilde{n}) =
(1\times\theta\times1)^* f_*u =\psi_*(i^{(3)}_*\tilde{\varphi}^*u)
\]
and hence $q_{12}^*\tilde{m} \otimes q_{23}^*\tilde{n}
=i^{(3)}_*\tilde{\varphi}^*u$. It follows that
\[
q_{13*}(q_{12}^*\tilde{m} \otimes q_{23}^*\tilde{n})
=q_{13*}i^{(3)}_*\tilde{\varphi}^*u
=i^{(2)}_*\overline{q}_{13*}\tilde{\varphi}^*u
=i^{(2)}_*\varphi^*(\overline{p}_{13*}u).
\]
By
\[
k_*(\overline{p}_{13*}u) =\tilde{p}_{13*}f_*u =\tilde{p}_{13*}
(\tilde{p}_{12}^*k_*m \otimes \tilde{p}_{23}^*k_*n)
\]
we conclude that the element of
$K^{G\times\BC^*}(Y^{(2)};Y_O^{(2)})$ corresponding to $m\bigstar
n\in K^{G\times\BC^*}(Z_O;Z_O)$ is given by
$q_{13*}(q_{12}^*\tilde{m} \otimes q_{23}^*\tilde{n})$.
Proposition \ref{prop:conv} is verified. This completes the proof
of Theorem \ref{thm:conv}.

\medskip

{\bf Acknowledgement:} The work started during TT's fruitful visit
to the Chinese Academy of Science in 2001,  and was completed
during NX's enjoyable visit to the Department of Mathematics,
Osaka City University in 2004 under the financial support of the
21-st century COE program ``Constitution of wide-angle
mathematical basis focused on knots". TT and NX would like to
thank the hospitality of those institutions.

\bibliographystyle{unsrt}

\begin{thebibliography}{99}
\bibitem{CG}
Chriss, N., Ginzburg, V.: Representation theory and complex
geometry, Birkh\"auser Boston, Inc., Boston, MA, 1997.
\bibitem{G}
Ginzburg, V.: Lagrangian construction of representations of Hecke
algebras, Adv. in Math. \textbf{63} (1987), 100--112.
\bibitem{G2}
Ginzburg, V.: Geometrical aspects of representation theory,
Proceedings of the International Congress of Mathematicians, Vol.
\textbf{1, 2} (Berkeley, Calif., 1986), 840--848, Amer. Math.
Soc., Providence, RI, 1987.
\bibitem{KL1}
Kazhdan, D., Lusztig, G.: Representations of Coxeter groups and
Hecke algebras,  Invent. Math. \textbf{53}  (1979), 165--184.
\bibitem{KL2}
Kazhdan, D., Lusztig, G.: Proof of the Deligne-Langlands
conjecture for Hecke algebras, Invent. Math. \textbf{87}  (1987),
153--215.
\bibitem{LC}
Lusztig, G.: Some problems in the representation theory of finite
Chevalley groups, The Santa Cruz Conference on Finite Groups
(Univ. California, Santa Cruz, Calif., 1979),  313--317, Proc.
Sympos. Pure Math., \textbf{37}, Amer. Math. Soc., Providence,
R.I., 1980.
\bibitem{LS}
Lusztig, G.: Some examples of square integrable representations of
semisimple $p$-adic groups, Trans. Amer. Math. Soc. \textbf{277}
(1983), 623--653.
\bibitem{L:cell}
Lusztig, G.: The two-sided cells of the affine Weyl group of type
$A\sb n$, Infinite-dimensional groups with applications (Berkeley,
Calif., 1984),
  275--283,
Math. Sci. Res. Inst. Publ., \textbf{4}, Springer, New York, 1985.
\bibitem{L:cell1}  Lusztig, G.:  Cells in affine Weyl groups, in
``Algebraic groups and related topics", Advanced Studies in Pure
Math., vol. \textbf{6}, Kinokuniya and North Holland, 1985, pp.
255-287.
\bibitem{L:cell2}  Lusztig, G.:  Cells in affine Weyl groups, II,
J. Alg. \textbf{109} (1987), 536-548.
\bibitem{L:cell4}
Lusztig, G.: Cells in affine Weyl groups. IV, J. Fac. Sci. Univ.
Tokyo Sect. IA Math. \textbf{36}  (1989), 297--328.
\bibitem{L1}
Lusztig, G.: Bases in equivariant $K$-theory, Represent. Theory
\textbf{2}  (1998), 298--369 (electronic).
\bibitem{L2}
Lusztig, G.: Bases in equivariant $K$-theory. II, Represent.
Theory \textbf{3}  (1999), 281--353 (electronic).
\bibitem{Shi}
Shi, J.-Y.: The Kazhdan-Lusztig cells in certain affine Weyl
groups, Lecture Notes in Mathematics, \textbf{1179}.
Springer-Verlag, Berlin, 1986.
\bibitem{Shi96}
Shi, J.-Y.: The partial order on two-sided cells of certain affine
Weyl groups, J. Algebra \textbf{179}  (1996),  607--621.
\bibitem{S}
Slodowy, P.: Simple singularities and simple algebraic groups,
Lecture Notes in Mathematics, \textbf{815}. Springer, Berlin,
1980.
\bibitem{T}
Tanisaki, T.: Hodge modules, equivariant $K$-theory and Hecke
algebras, Publ. Res. Inst. Math. Sci. \textbf{23}  (1987),
841--879.
\bibitem{T2}
Tanisaki, T.: Representations of semisimple Lie groups and
$D$-modules, Sugaku expositions \textbf{4}  (1991), 43--61.
\bibitem{Th}
Thomason, R. W.: Equivariant algebraic versus topological
$K$-homology Atiyah-Segal style, Duke Math. J.  \textbf{56}
(1988), 589--636.
\bibitem{Xi}
Xi, N.: Representations of affine Hecke algebras, Lecture Notes in
Mathematics, \textbf{1587}. Springer-Verlag, Berlin, 1994.
\bibitem{X:b}  Xi, N.:  The based ring of two-sided cells of
affine Weyl groups of type $\tilde A_{n-1},$ Mem. of AMS, Vol.
157, No.749, 2002.
\end{thebibliography}

\end{document}